\documentclass[12pt,leqno]{amsart}
\usepackage{latexsym, amssymb, amsmath, amsfonts}

\begin{document}
\newtheorem{theorem}{Theorem}[section]
\newtheorem{proposition}[theorem]{Proposition}
\newtheorem{fact}[equation]{}
\newtheorem{corollary}[theorem]{Corollary}
\newtheorem{definition}[theorem]{Definition}
\newtheorem{remark}[theorem]{Remark}
\newtheorem{remarks}[theorem]{Remarks}
\newtheorem{lemma}[theorem]{Lemma}
\newtheorem{maintheorem}[theorem]{Main Theorem}
\newenvironment{pf}{\proof[\proofname]}{\endproof}

\newtheorem{Thm}{Theorem}[section]
\newtheorem{Cor}[Thm]{Corollary}
\newtheorem{Rem}[Thm]{Remark}
\newcommand{\B}[1]{\mathbb#1}
\newcommand{\cal}[1]{\mathcal{#1}}
\newcommand{\C}[1]{\cal#1}
\newcommand{\isom}{\overset {\thicksim}{\to}}
\newcommand{\oom}[2]{\text{$\Omega^{#1}_{#2}$}}
\newcommand{\Om}{\oom{d}{K_w}}
\newcommand{\wt}{\widetilde}
\newcommand{\GM}{\Gamma}
\newcommand{\gm}{\gamma}
\newcommand{\dt}{\delta}
\newcommand{\Dt}{\Delta}
\newcommand{\bs}{\backslash}
\newcommand{\m}{^{\times}}
\newcommand{\Fr}{\text{Fr}}


\newcommand{\set}[1]{\{\,#1\,\}}
\newcommand{\pair}[2]{\langle #1,#2 \rangle}
\newcommand{\matr}[4]{\left[
\begin{array}{rr} #1 & #2 \\ #3 & #4
                    \end{array} \right]}
\newcommand{\smatr}[4]{\left[
\begin{array}{rr}
{\scriptstyle #1} & {\scriptstyle #2}\\
{\scriptstyle #3} & {\scriptstyle #4}
\end{array} \right]}
\newcommand{\ov}{\overline}
\newcommand{\ra}{\rightarrow}
\newcommand{\la}{\leftarrow}
\newcommand{\lra}{\longrightarrow}
\newcommand{\lla}{\longleftarrow}
\newcommand{\ua}{\uparrow}
\newcommand{\da}{\downarrow}
\newcommand{\sra}[1]{\stackrel{#1}{\ra}}
\newcommand{\sla}[1]{\stackrel{#1}{\la}}
\newcommand{\nr}{{{\rm nr}}}
\newcommand{\rint}{{{\rm int}}}
\newcommand{\rf}{{\rm f}}
\newcommand{\rref}[1]{{\rm (\ref{#1})}}
\newcommand{\Aut}{\operatorname{Aut}}
\newcommand{\Cl}{\operatorname{Cl}}
\newcommand{\diag}{\operatorname{diag}}
\newcommand{\Disc}{\operatorname{Disc}}
\newcommand{\dist}{\operatorname{dist}}
\newcommand{\Emb}{\operatorname{Emb}}
\newcommand{\End}{\operatorname{End}}
\newcommand{\Ed}{\operatorname{Ed}}
\newcommand{\Gal}{\operatorname{Gal}}
\newcommand{\Gr}{\operatorname{Gr}}
\newcommand{\Hom}{\operatorname{Hom}}
\newcommand{\Id}{\operatorname{Id}}
\newcommand{\inv}{\operatorname{inv}}
\newcommand{\Jac}{\operatorname{Jac}}
\newcommand{\Ker}{\operatorname{Ker}}
\newcommand{\Map}{\operatorname{Map}}
\newcommand{\Mat}{\operatorname{Mat}}
\newcommand{\Nm}{\operatorname{Nm}}
\newcommand{\res}{\operatorname{res}}
\newcommand{\Res}{\operatorname{Res}}
\newcommand{\Spec}{\operatorname{Spec}}
\newcommand{\Stab}{\operatorname{Stab}}
\newcommand{\Tr}{\operatorname{Tr}}
\newcommand{\val}{\operatorname{val}}
\newcommand{\Var}{\operatorname{Var}}
\newcommand{\Ver}{\operatorname{Ver}}

\newcommand{\SL}{\operatorname{SL}}
\newcommand{\SU}{\operatorname{SU}}
\newcommand{\GL}{\operatorname{GL}}
\newcommand{\PGL}{\operatorname{PGL}}
\newcommand{\GU}{{\bf{GU}}\,}

\newcommand{\cA}{{\mathcal A}}
\newcommand{\cE}{{\mathcal E}}
\newcommand{\cF}{{\mathcal F}}
\newcommand{\cG}{{\mathcal G}}
\newcommand{\cH}{{\mathcal H}}
\newcommand{\cK}{{\mathcal K}}
\newcommand{\cL}{{\mathcal L}}
\newcommand{\cM}{{\mathcal M}}
\newcommand{\cO}{{\mathcal O}}
\newcommand{\cP}{{\mathcal P}}
\newcommand{\cQ}{{\mathcal Q}}
\newcommand{\cS}{{\mathcal S}}
\newcommand{\cX}{{\mathcal X}}
\newcommand{\cY}{{\mathcal Y}}
\newcommand{\bbA}{{\Bbb A}}
\newcommand{\bbC}{{\Bbb C}}
\newcommand{\bbF}{{\Bbb F}}
\newcommand{\bbQ}{{\Bbb Q}}
\newcommand{\bbO}{{\Bbb O}}
\newcommand{\bbR}{{\Bbb R}}
\newcommand{\bbZ}{{\Bbb Z}}
\newcommand{\bbG}{{\Bbb G}}
\newcommand{\bbH}{{\Bbb H}}
\newcommand{\bbP}{{\Bbb P}}
\newcommand{\bA}{{\bf A}}
\newcommand{\bB}{{\bf B}}
\newcommand{\bF}{{\bf F}}
\newcommand{\bG}{{\bf G}}
\newcommand{\bH}{{\bf H}}
\newcommand{\bM}{{\bf M}}
\newcommand{\bZ}{{\bf Z}}
\newcommand{\lam}{\lambda}

\newcommand\qi{\hat{\imath}}
\newcommand\qj{\hat{\jmath}}
\newcommand\qk{\hat{k}}
\renewcommand{\ast}{\times}
\newcommand{\Clplus}{{\rm Cl}^{+}}
\newcommand{\Dint}{D^\rint}
\newcommand{\Bint}{B^{\rint}}
\newcommand{\alint}{\alpha^{\rint}}
\newcommand{\bGint}{\bG^{\rint}}
\newcommand{\nuint}{\nu^\rint}
\newcommand{\Nmint}{\Nm^\rint}
\newcommand{\Shim}{X^\rint}
\newcommand{\ShimS}
   {S\backslash \Shim}
\newcommand{\ShimCS}
   {\ShimS_{\bbC}}
\newcommand{\ShimSan}
   {(\ShimCS)^{\rm an}}
\newcommand{\Shimp}{X^{\rint,p}}
\newcommand{\Ttil}{{\tilde{T}}}
\newcommand{\Gm}{\bbG_m}
\newcommand{\Frob}{{\rm Frob}}
\newcommand{\cOstFP}{\cO_{F,\cP}^\times}
\newcommand{\Minfp}{M_\infty^+}
\newcommand{\bcdot}{{\bf \cdot}}
\newcommand{\Omhat}{\hat{\Omega}
            \vphantom{\Omega}}
\newcommand{\Edo}{\Ed^{\rm o}}
\newcommand{\invlim}
     {\displaystyle{\lim_{\leftarrow}}}
\newcommand{\XLO}{X_{L,0}}
\newcommand{\Kgal}{K_{\rm Gal}}

\title{Local Points on $P$-adically Uniformized
Shimura Varieties}
\author{Bruce W. Jordan}
\address{Department of Mathematics, Box G-0930\\
Baruch College, CUNY\\
17 Lexington Avenue\\
New York, NY  10010  USA}
\thanks{The first author was partially supported
by grants from the NSF and PSC-CUNY}
\author{Ron Livn\'{e}}
\address{Mathematics Institute\\
The Hebrew University of Jerusalem\\
Givat-Ram, Jerusalem  91904\\
Israel}
\email{rlivne@sunset.math.huji.ac.il}
\thanks{The first two authors were partially
supported by a joint Binational Israel-USA
Foundation grant}
\author{Yakov Varshavsky}
\address{Mathematics Institute\\
The Hebrew University of Jerusalem\\
Givat-Ram, Jerusalem  91904\\
Israel}
\email{vyakov@math.huji.ac.il}
\date{March 7, 2002}
\keywords{Shimura varieties, local points,
$p$\/-adic uniformization, even jacobians}
\subjclass{Primary: ; Secondary: }
\begin{abstract}
Using the $p$\/-adic uniformization of Shimura
varieties we determine, for some of them, over which
local fields they have rational points. Using this
we show in some new curve cases that the jacobians
are even in the sense of \cite{PS}.
\end{abstract}

\maketitle

\section{Introduction and notation}
The real points on certain Shimura varieties
were determined by Shimura in \cite{Shi}.
In \cite{JL1} the case of
$p$\/-adic points was treated for Shimura
curves associated to maximal orders in
indefinite rational quaternion division
algebras. The case of good reduction turned
out to reduce to the trace formula via
Hensel's lemma. The case of bad reduction was
handled through the result of \v{C}erednik and
Drinfel'd, which assured that these curves
admit a $p$\/-adic uniformization.

A general result of Varshavsky, Rapoport and Zink
\cite{Va1, Va2} gives the known cases when Shimura
varieties admit a $p$\/-adic uniformization. Of
special interest in this context is the case of
curves, because of the recent results of \cite{PS}
and \cite{JL3}. In this work we answer the question
of existence of local points for some of these
varieties. Using this we show in some new curve
cases that the jacobians are even in the sense of
\cite{PS}.

We now set some notation. Let $F$ be a totally real
number field, $g=[F:\bbQ]$. Denote by
$\Sigma_{\infty}=\{\infty_1,\ldots,\infty_g\}$ the
set of real embeddings of $F$ and by $\Sigma_\rf$
the set of finite places of $F$. Set $\Sigma =
\Sigma_{\infty}\cup \Sigma_\rf$. Let $\bbA = \bbA_F$,
$\bbA_\rf$, $\bbA_\rf^{\Xi}$, $\bbO_\rf^{\Xi}$, $\bbA^\times$,
$\bbA_\rf^{\times}$, and $\bbO_\rf^{\times}$ denote the ad\`eles of $F$,
the finite ad\`eles of $F$, the finite ad\`{e}les
without the $v$ components for places $v$ in a
finite set of places $\Xi\subseteq\Sigma_{f}$, the integral finite ad\`{e}les
without the $v$ components for places $v$ in a finite set of places $\Xi\subseteq\Sigma_{f}$, the
id\`{e}les of $F$, the finite id\`eles of $F$, and the integral finite id\`eles of $F$,
respectively.  For an algebraic group $\bG/F$ we
abbreviate $G=\bG(F)$, $G_{v}=\bG(F_{v})$,
$G_{\infty}=\prod_{i=1}^{g}\bG(F_{\infty_i})$,
$G_\rf = \bG(\bbA_\rf)$, and
$G_\rf^{\Xi}=\bG(\bbA_\rf^{\Xi})$.
 We view $G$ as
contained in $\bG(\bbA)$ and also in each
$G_\rf^{v_1,\ldots , v_r}$
and in each $G_{v}$ via the natural projection.

\section{$P$\/-adic uniformization of Shimura
varieties}
\label{padic}
 We will consider two types of Shimura varieties
(cf.\ \cite{De2}, \cite{Va2}):\\[.05in]
{\bf Case 1}.\hspace*{1em}Let $\Bint/F$ be a
quaternion division algebra, Let $\bGint /F$ be the
algebraic group associated to its multiplicative
group, and let $\nuint:\bGint\ra\Gm$ be the
$F$-morphism induced from the reduced norm from
$\Bint$ to $F$.  Assume that
$\Bint\otimes_{F,\infty_i}\bbR$ is indefinite for
$1\leq i\leq r$ and definite for $r+1\leq i\leq g$,
with $r\geq 1$. Fix identifications of
$B_{\infty_i}$ with $\Mat_{2\times2}(\bbR )$ for
$1\leq i\leq r$ and with the Hamilton quaternions
$\bbH$ for $r+1\leq i\leq g$.  Then
$B\otimes_{\bbQ}\bbR\cong
\Mat_{2\times2}(\bbR )^{r}\times \bbH^{g-r}$.
Let $h:\Res_{\bbC/\bbR}\Gm(\bbR)\rightarrow
\Res_{F/\bbQ}\bGint (\bbR) $ be the Hodge type
\[
h(z)= (m(z), m(z),\ldots, m(z),1,1,\ldots, 1) \in
\GL_{2}(\bbR)^{r}\times (\bbH^\times)^{g-r} =
\Res_{F/\bbQ}\bGint (\bbR)\ ,
\]
where for
$z=x+\sqrt{-1}y \in\bbC^{\times} =
             \Res_{\bbC/\bbR}\Gm(\bbR)$
we set $m(z) = \smatr{x}{y}{-y}{x}^{-1}$. \\[.05in]
{\bf Case 2}.\hspace*{1em}Here let $K$ be a CM
extension of $F$, let $(x\mapsto \ov{x}):K\ra K$
be the conjugation over $F$, and let $d\geq 2$
be an integer. Let $\Dint$ be a division algebra,
central and of dimension $d^2$ over $K$,
with an involution of the second kind $\alint$.
Let $\bGint={\bf GU}(\Dint,\alint )$ be the
associated group of unitary similitudes and let
$\nuint:\bGint\ra\Gm$ be the $F$-morphism
induced from the factor of similitudes. Also  let
$\Nmint:\bGint\ra\Res_{K/F}\Gm$ be the morphism
induced from the norm map from $\Dint$ to $K$. Choose
extensions $\infty_{1},\ldots,\infty_{g}$ to
embeddings, denoted by the same letter,
$\infty_i:K\rightarrow \bbC$, and suppose that
\[
G_{\infty_i} \cong \left\{
\begin{array}{l}
\GU_{d-1,1}(\bbR ) \quad\mbox{for}
              \quad 1\leq i\leq r\ ,\\
\GU_{d}(\bbR ) \quad\mbox{for}\quad r+1\leq i \leq g
\end{array}
\right. .
\]
Define  $h:\Res_{\bbC/\bbR}\Gm (\bbR )\rightarrow
\Res_{F/\bbQ}\bGint(\bbR) \cong
\GU_{d-1}(\bbR )^{r}\times \GU _{d}(\bbR )^{g-r}$
by
\[
h(z) = (m(z), \ldots , m(z), 1, \ldots, 1)\in
\GU_{d-1,1}(\bbR)^r\times \GU_{d}(\bbR )^{g-r}\ ,
\]
where $m(z)=\diag(1,\ldots, 1,z/\ov{z})^{-1}$.

To handle cases 1 and 2 simultaneously, we formally
set $K=F$ and $d=2$ in Case $1$.

In both cases, let $\Kgal$ be
the Galois closure of $\infty_{1}(K)$
in $\bbC$, and let $H$ be the subgroup of
$\Gal (\Kgal/\bbQ)$ preserving the set
$\{\infty_{1},\ldots, \infty_{r}\}$.  The
reflex field is defined by $E=\Kgal^{H}$.

Let $X_{\infty}$ be the conjugacy class of
$h$. Our identifications allow us to view
$X_{\infty}$ as $(\bbC\smallsetminus\bbR)^{r}$ in the case $d=2$, 
and as $(B^{d-1})^{r}$ in the case $d\geq 3$, where $B^{d-1}$ is the unit ball
in $\bbC^{d-1}$. For a compact
open subgroup $S\subset G_\rf$, the
corresponding Shimura variety is given complex
analytically by
\[
\wt{X^\rint_S}= \bigl(X_{\infty}\times
(S\backslash G^{\rint}_\rf)\bigr)
\bigm/ G^{\rm int}\ ,
\]
where  $G^{\rm int}$ acts on the product $X_{\infty}\times
(S\backslash G^{\rint}_{\rf})$ by the rule $(x,g)\gm:=(\gm^{-1}(x),g\gm)$.
This variety is denoted by $_{S}M_{\bbC}(H,X_{\infty})$ in
\cite{De1}.

The complex manifold $\wt{X^\rint_S}$ has a canonical
algebraization $X^\rint_S$ to a complex
projective variety.  In fact each $X^\rint_S$
admits a canonical model over the reflex field $E$ and $X^\rint_S$'s form a projective system.
The main result of \cite{Va2}
is that $X^\rint=\invlim_{S} X^\rint_S$ admits a
$p$-adic uniformization, under some conditions
which we will now specify. Let $v$ be
a finite place of $E$ of residue characteristic
$p$ and choose  a prime $\tilde{v}$ of $\Kgal$
above $v$.  Fix an embedding of the completion
$\hat{K}_{{\rm Gal},\tilde{v}}$ of
$\Kgal$ into $\bbC$. The embeddings
$\infty_1 , \ldots ,\infty_r$ then induce
places $w_1,\ldots,w_r$ of $K$ above $p$.
Moreover $\set{w_1,\ldots , w_r}$ depends only
on $v$, because $\Gal (\Kgal /E)$ preserves
$\set{\infty_{1},\ldots , \infty_{r}}$, and
$E_v$ is the compositum of
$K_{w_1},\ldots, K_{w_r}$ in
$K_{{\rm Gal},\tilde{v}}$ as in
\cite[Lemma 2.6]{Va2}. Let $v_i$ be the
restriction of $w_i$ to $F$. Then
$X_{v}^\rint=X^\rint\otimes_{E}E_{v}$ has a
$p$-adic uniformization (\cite[Theorems 5.3 and
2.13]{Va2}) provided we have
\begin{quote}
1.\hspace*{1em}The $v_i$'s are distinct and
split in $K/F$.\\
2.\hspace*{1em}The Brauer invariant
$\inv_{w_i\,}D^\rint$ is $1/d$ for each
$1\leq i\leq r$.
\end{quote}
We shall now describe a special case of the
result more precisely. This requires more
notation.

In Case 1, let $B$ be a quaternion algebra over
$F$ split at $v_1$,\ldots , $v_r$; ramified at
all infinite places; and with a fixed
isomorphism
\[
B^{\rm int}\otimes \bbA_{\rf}^{v_1,\ldots, v_r}
\cong
B\otimes \bbA_\rf^{v_1,\ldots, v_r}\ .
\]
Let $\bG$ be the algebraic group over $F$
corresponding to $B^{\times}$, and let
$\nu:\bG\ra\Gm$ be the $F$-morphism induced
from the reduced norm from $B$ to $F$. Set
$\cG=\cG'=G_\rf^{v_1,\ldots , v_r}$.

In Case 2, let $D$ be a central simple
$K$-algebra with an involution $\alpha$ of the
second kind.  We suppose that $D$ is split at
all places over $v_1,\ldots ,v_r$ (which are
$w_1, {\ov{w}}_1,\ldots ,w_r, {\ov{w}}_r$),
and that we are given an identification of
$(D, \alpha) \otimes_F F_v$ with
$(D^{\rm int}, \alpha^{\rm int})\otimes_F F_v$
at the other finite places $v$ of $F$. Moreover we
assume that $\alpha$ is definite at all the
infinite places of $K$. Define an algebraic
group $\bG$ over $F$ by $\bG = \GU(D,\alpha)$,
and let $\nu:\bG\ra\Gm$ be the $F$-morphism
induced from the factor of similitudes. Also
let $\Nm:\bG\ra\Res_{K/F}\Gm$ be the
morphism induced from the norm map from $D$ to
$K$. Set $\cG = G_\rf^{v_1,\ldots,v_r}$. The
decompositions
\[
D\otimes_F F_{v_i}  \simeq  \Mat_d(K_{w_i})\oplus
\Mat_d(K_{\ov{w_i}})\quad\mbox{and}\quad
D^\rint \otimes_F F_v  \cong
D^\rint_{w_i}\oplus D^\rint_{\ov{w_i}}\,,
\]
together with the $F$-rational maps
$\nu:\bG\rightarrow \Gm$ and
$\nuint:\bG^{\rm int} \rightarrow \Gm$ induce
decompositions
\begin{equation}
\label{decG}
G_{v_i}\cong\GL_{d}(K_{w_i})\times
F_{v_i}^{\times} \qquad \mbox{and}\qquad
G_{v_i}^{\rm int}\cong (D^\rint_{w_i})^{\times}\times
F_{v_i}^{\times}\ .
\end{equation}
More precisely, the conjugation of $K$ over
$F$ induces an isomorphism
$K_{w_i}\simeq K_{\ov{w_i}}$ and $\alpha$
(respectively $\alpha^\rint$) induces a
compatible isomorphism of algebras
$D_{\ov{w_i}}\simeq D_{w_i}^{{\rm opp}}$
(respectively
$D^\rint_{\ov{w_i}}\simeq
             (D_{w_i}^\rint)^{{\rm opp}}$).
Viewing all the above isomorphisms as
identifications, the
decompositions~\rref{decG} above become
\begin{equation}
\label{decbis}
\quad G_{v_i} = \{(g,\lambda g)|
g\in\GL_d(K_{w_i}),\;\lambda\in
F^\times_{v_i}\} \quad \mbox{and}\quad
G^\rint_{v_i} = \{(g,\lambda g)| g\in (D^\rint_{w_i})\m,
          \;\lambda\in F^\times_{v_i}\}\ .
\end{equation}
We can then view
$\cG' =  \prod_i F_{v_i}^{\times} \times \cG$
as a {\em closed} subgroup of both $G_\rf$ and
$G_\rf^{{\rm int}}$.

In both cases, fix
for each $i$ a central division algebra
$\wt{D}_{w_i}$  over $K_{w_i}$ of invariant $1/d$
($=1/2$ in Case 1). We then have that
\[
G_\rf \simeq \prod_{i=1}^r \GL_d(K_{w_i})
\times \cG' \qquad\text{and}\qquad
G_\rf^{{\rm int}} \simeq
\prod_{i=1}^r \wt{D}_{w_i}^\ast \times \cG'\ .
\]

The ring of integers $\cO_{\wt{D}_{w_i}}$  is
normalized by $\wt{D}_{w_i}^{\times}$, and
$\wt{D}_{w_i}^{\times}/\cO_{\wt{D}_{w_i}}^\times$ is
identified with $\bbZ$ via the valuation of
the norm.  Hence the group $G_\rf^{\rm int}$
acts on $\Shimp:=(\prod_{i=1}^r
   \cO_{\wt{D}_{w_i}}^{\times}) \backslash\Shim$ through its quotient
\[ G_\rf^{{\rm int}}\bigm/
     \prod_{i=1}^r \cO_{\wt{D}_{w_i}}^\times
       \cong\cG'  \times \bbZ^{r}\ . \]

For a finite extension $L$ of $\bbQ_p$, let
$\Omega_L = \Omega_L^d$ be Drinfeld's symmetric
space.  Recall that the analytic space (\cite{Ber}) $\Omega_L$ is
obtained by removing all rational hyperplanes from
$\bbP^{d-1}_L$. It is preserved by the action of
$\PGL_d(L)$. Consider an analytic space $\Omega_L^\nr=\Omega_L\hat{\otimes}_L
L^\nr$ over $L$, where $L^\nr$ is the completed maximal
unramified extension of $L$.  We let $g\in\GL_d(L)$
act on $\Omega_L^\nr$ via the natural (left) action on
$\Omega_L$ and the action of $\Frob_L^{\val(\det
g)}$ on $L^\nr$. We also let $n\in \bbZ$ act on
$\Omega_L^\nr$ through the action of $\Frob_L^{-n}$
on $L^\nr$. This gives an  $L$-rational action of $\GL_d(L)\times
\bbZ$ on $\Omega_L^\nr$. In our situation, set
\[
\Omega = \prod_{i=1}^r\Bigl(\Omega_{K_{w_i}}^\nr
\hat{\otimes}_{K_{w_i}}E_v\Bigr)\ .
\]
This is an $E_v$-analytic space with an $E_v$-rational
action of the group
$\prod_{i=1}^{r}\GL_d(K_{w_i})\times \bbZ^r$.

Let $G$ act on the $i^{\rm th}$ factor of $\Omega$
via the embedding
$G(F)\hookrightarrow \GL_d(K_{w_i})$, and
on $\cG'$ through the natural embedding
(in Case 2, let $G(F)$ act on each factor
$F_{v_i}^{\times}$ of $\cG'$ via the
similitude factor $\nu$).

For a compact open subgroup $S\subset \cG'$,
the $E_v$\/-analytic space
$ (\Omega\times (S\backslash \cG'))/G(F)$
algebraizes canonically to a scheme
$X_S$ over $E_v$. The inverse limit
of $X_S$ over such $S$'s is a scheme
$X$ over $E_v$ with an $E_v$-rational action
by $\cG' \times \bbZ^r$.

A special case of the main result of
\cite{Va2} is the following
\begin{theorem}
\label{int}
Under conditions {\rm 1} and {\rm 2} above,
there exists a $\cG' \times \bbZ^r$-equivariant,
$E_v$-rational isomorphism
\[   \Shimp\otimes_E E_v\cong X \,.    \]
\end{theorem}

\section{The connected component}
\label{concomp}
 A description of the set of
connected components of a Shimura variety, with the
Galois action, is known in general (see
\cite[2.1.3.1]{De1}). We need an explicit form of
this description in our case(s). For this it is
simpler to use the $p$\/-adic uniformization,
essentially repeating the standard argument there.

Define a torus over $F$ by
$\bM = \bG/\bG_{{\rm der}}$. In Case 1, $\nu$
induces an isomorphism $\bM\simeq\Gm/F$; in
Case 2, $\Nm\times\nu$ induces an isomorphism
of $\bM$ with the torus over $F$ defined by
\[
\{(x,\lambda)\in \Res_{K/F}\Gm \times
            \Gm | x\ov{x} = \lambda^d\}\,,
\]
(Notice that $\bM$ is also isomorphic to
$\bGint/\bGint_{{\rm der}}$ through $\nuint$
and $\Nmint$.)

In Case 1, $M_{v_i}\simeq F^\times_{v_i}$; in
Case 2, we have
$K\otimes_F F_{v_i}\cong
                 K_{w_i}\oplus K_{\ov{w_i}}$,
and therefore we can and will make the
identification
\begin{equation}
\label{decM}
M_{v_i} = K^\times_{w_i}\times F^\times_{v_i}\,.
\end{equation}
Using the decomposition~\rref{decG}, the
quotient map $G_{v_i}\ra M_{v_i}$ is then
induced by $\Nm\times \nu$.

Since the derived group $\bG_{\rm der}$
is a form of $\SL_2$ in the first
case and ${\bf SU}_d$ in the second case,
it is simply connected.  Using Hasse
principle, combined with vanishing of cohomology of $p$-adic groups, and a
local archimedean calculation, we get that the image of
$G$ in $M$ is
\begin{equation}
\label{hasse}
M^+ = \{x\in M| \nu(x)\in F^\times
\text{ is positive at all the infinite
places of $F$ }\}.
\end{equation}

For a compact open subgroup $S\subset \cG'$,
let $T$ be the image of $S$ in $M_\rf$, and
set
${\wt{T}} = \Big( \prod_{i=1}^r
       \cO^\times_{K_{w_i}} \Big) \times T$.
Then ${\wt{T}}$ is an open subgroup of $M_\rf$.
Let $\Minfp$ be the connected
component of $\bM(F \otimes_\bbQ \bbR)$. Set
$\Cl^+_{\wt{T}} =
      (\Minfp{\wt{T}})\backslash \bM(\bbA)/ M$.
Let $f_i = f(E_v/K_{w_i})$ be the degree of the
extension of residue fields. Let $\pi_{w_i}$ be a
uniformizer of $K_{w_i}$. Using the
decomposition~\rref{decM}, we can view $\pi_{w_i}$
as an element of $M_{v_i}$. The image $\varphi_+$
of $\prod_{i=1}^r \pi_{w_i}^{f_i}$ in
$\Cl^+_{\wt{T}}$ is then independent of the choices
of the $\pi_{w_i}$\/'s.

For a local field $L$, denote  by $L^{(r)}$
the extension of degree $r$ of $L$ in $L^\nr$.
We will now prove the following:
\begin{theorem}
\label{fieldofdef}
The map $\nu$ in Case {\rm 1} and
$\Nm\times\nu$ in
Case {\rm 2} induce an isomorphism
$\pi_0(S\backslash \Shimp)
                      \sra{\sim} \Cl^+_{\wt{T}}$.
Let $k_+$ be the order of $\varphi_+$ in
$\Cl^+_{\wt{T}}$. Then the field of definition of
each geometric connected component of
$S\backslash \Shimp$\/ is $E_v^{(k_+)}$.
\end{theorem}
\begin{pf}
Fix an algebraic closure $\ov{E_v}$ of $E_v$.
Since each
$\Omega_{K_{w_i}} \hat{\otimes}_{K_{w_i}}E_v$
is geometrically connected, the set of
geometrically connected components of
$\Omega^\nr_{K_{w_i}}
               \hat{\otimes}_{K_{w_i}}E_v$
is in bijection with the set $\Emb_i$ of continuous 
$K_{w_i}$\/-homomorphisms of $K_{w_i}^\nr$
into the completion of $\ov{E_v}$. Each $\Emb_i$ carries a
natural $\Gal(K_{w_i}^\nr/K_{w_i}) \times
\Gal(\ov{E_v}/E_v)$ action. The two Galois
actions are obviously related: the
$\Gal(\ov{E_v}/E_v)$\/-action factors through
$\Gal(K_{w_i}^\nr E_v/E_v) = \Gal(E_v^\nr/E_v)$, and
\begin{equation}
\label{Frobs}
\Frob_{E_v} \in \Gal(E_v^\nr/E_v)\; \text{ acts
like } \;\Frob_{K_{w_i}}^{-f_i} \in
\Gal(K_{w_i}^\nr/K_{w_i})\ .
\end{equation}
The $\Gal(K_{w_i}^\nr/K_{w_i})$\/-action induces
an action of $G\subset \GL(2,K_{w_i})$ through
\begin{equation}
\label{GandGal}
\GL(2,K_{w_i})\ni g\mapsto
     \Frob_{K_{w_i}}^ {\val_{w_i}(\det g)}\ .
\end{equation}
The set of connected components
$\pi_0(S\backslash X^{\rint,p})$ is therefore
canonically the separated quotient
$(\prod_{i=1}^r\Emb _i\times
                        (S\backslash\cG'))/G$.
By the strong approximation theorem,
$G_{\rm der,\rf}^{ v_1, \ldots ,v_r}
            \subseteq G_{\rm der}S$, and
$G_{\rm der}$ acts trivially on
$\prod_{i=1}^{r}\Emb_i$. Therefore
\[
\pi_{0}(S\backslash X^{{\rm int},p}) =
 \Bigl( \bigl( \prod_{i=1}^{r} \Emb_i \bigr)
   \times (S\backslash \cG ')\Bigr)/G
    \cong \bigl( \prod_{i=1}^{r}\Emb_i\times
        (T\backslash M')\bigr)\bigm/ M^+\ ,
\]
where $M'$ is the image of $\cG'$ in $M_\rf$.
Explicitly, $M' = M_\rf^{ v_1,\ldots , v_r}$
in Case 1 and
$M' = M_\rf^{ v_1,\ldots , v_r}\times
\prod_{i=1}^r F^\times_{v_i}$ in Case 2. Now
fix some embedding of each $K_{w_i}^\nr$ in
$\ov{E_v}$. Then $\Emb_i$ gets identified with
$\Gal (K_{w_i}^\nr/K_{w_i})\cong\hat{\bbZ}$,
into which
$K_{w_i}^{\times}/\cO_{w_i}^{\times}$ embeds
as a dense subset preserved by the $G$-action
(through the $w_i$\/-valuation of the norm).
Using \rref{decM}, we get an embedding of
$M_v/\cO^\times_{K_{w_i}}
    \cong K^\times_{w_i}/\cO^\times_{K_{w_i}}\times F^\times_{v_i}$
into $\Emb_i\times F^\times_{v_i}$ as a dense
subset, compatible with the $M^+$\/-action
(through the $K_{w_i}$\/-component). Hence
\[
\pi_0(S\backslash \Shimp) \cong
\bigl(\prod_{i=1}^r
     (K^\times_{w_i}/\cO^\times_{w_i})
 \times (T\backslash M')\bigr) \bigm/ M^+ =
{\wt{T}}\backslash M_\rf/M^+
\cong {\wt{T}} \Minfp\backslash \bM(\bbA)/M\ ,
\]
using the weak approximation at the infinite
places to get the last isomorphism. Of course
this agrees with Deligne's general description
cited above, but now we see that
$\Gal (\ov{E_v}/E_v)$ acts on
$\pi_0(S\backslash \Shimp)$ through its image
in
$\prod_{i=1}^{r}\Gal(E_{v}K_{w_i}^\nr/E_{v})$,
and this image already factors through
$\Gal(E_v^\nr/E_v)$. By \rref{Frobs} and
\rref{GandGal}, $\Frob_{E_v}$ acts on
$\Cl^+_T$ as the id\`{e}le whose $v^{\text{th}}$
component is $\pi_{w_i}^{f_i}$ if $w=v_i$
and $1$ otherwise. We also see that a power
$\Frob_{E_{v}}^l$ acts trivially on one
component (equivalently, on all components)
if and only if $\varphi_+^l=1$ in $\Cl^+_T$.
This concludes the proof of the Theorem.
\end{pf}
\begin{remark}{\rm When $r=1$ we have
$E_v=K_{w_1}$, so that $\varphi_+ = \pi_{w_1}$.}
\end{remark}

\section{Local points of twisted Mumford quotients}
\label{locpoints}
 In this section only we change
our notation, and let $F$ be a finite extension
field of $\bbQ_p$. Let $d\geq 2$ be an integer, and
let $\GM\subset\GL_d(F)$ be a subgroup. Set
$Z\GM:=\GM\cap Z(\GL_d(F))$, and for every subgroup
$\Dt\subset\GM$ containing $Z\GM$, we will write
$P\Dt$ instead of $\Dt/Z\GM$. Assume that:

A) the closure of $\GM$ has a finite covolume in
$\GL_d(F)$;

B) $P\GM$ is a (cocompact) lattice in $\PGL_d(F)$.

C) (for convenience) $\det\GM\subset Z\GM$.

Let $\val(\det\GM)=k_+\B{Z}$ and
$\val(Z\GM)=k\B{Z}$. Then by C), $k|k_+|dk$. Let
$$\GM':=\{\gm\in\GM| dk \text{ divides
}\val(\det(\gm))\}.$$ Then $\GM'$ is a normal
subgroup of $\GM$ containing $Z\GM$, and $\GM/\GM'$ is cyclic of
order $dk/k_+$. For every subgroup $\Dt\subset\GM$
we denote by $\Dt'$ the intersection of $\Dt$ with
$\GM'$.

Let $X: = \GM\bs (\oom{d}{F}\hat{\otimes}_{F}{F^{nr}})$, and let $X':=\GM'\bs \oom{d}{F}$.
Then $X$ and $X'$ are projective and geometrically
connected varieties over $F^{(k_+)}$ and $F$
respectively, and
$X = (\GM/\GM')\bs (X'\otimes_F F^{(dk)})$. Thus $X$ is a {\em {\rm
(}Frobenius{\rm )} twist} of the {\em Mumford
uniformized} variety $X'$.

Let $L$ be any finite extension of $F^{(k_+)}$, and
let $e=e(L/F^{(k_+)})$ and $f=f(L/F^{(k_+)})$ be the
ramification and the inertia degrees of $L$ over
$F^{(k_+)}$ respectively. Notice that $[L:F] =
efk_+$.

\begin{Thm} \label{T:Main}
a) If $X(L)\neq\emptyset$, then there exists a
subgroup $\Dt\subset\GM$ containing $Z\GM$ such
that:

(i) the group $P\Dt$ is finite;

(ii) $\val(\det(\Dt))=dk\B{Z}+fk_+\B{Z}$;

(iii) $d|fek_+m$, where $m$ is the order of $P\Dt'$.

b) The converse of a) holds if we assume in
addition either that $d$ is prime or that
$d|efk_+$. In particular, the converse of a) holds
if $P\GM'$ is torsion-free.
\end {Thm}
\begin{pf}
Let $d'$ be the least integer $\geq 1$ for which
$L^{(d')}$ contains $F^{(dk)}$. Then $d'$ is the smallest positive integer such that
$dk$ divides $fk_+d'$, hence $d' = d/\gcd(d, fk_+/k)$. Observe that conditions (ii) and
(iii) of the theorem can be restated as (ii$'$)
$\val(\det(\Dt))=(dk/d')\B{Z}$ and (iii$'$) $d'| ekm$
respectively. Then $F^{(dk)} \otimes_{F^{(k_+)}} L
\simeq \oplus_{\Hom_{F^{(k_+)}}(F^{(kd)},L^{(d')})} L^{(d')}$. This
yields the following description of $X_L$\/:
\[X_L= X \otimes_{F^{(k_+)}} L \simeq
(\GM/\GM')\bs\left( \sqcup_{\Hom_{F^{(k_+)}}(F^{(kd)},L^{(d')})} X' \otimes_F
L^{(d')} \right)\simeq (\GM_1/\GM')\bs (X' \otimes_F
L^{(d')}),\] where $\GM_1$ is the
stabilizer in $\GM$ of any (and hence of each)
connected component of $X' \otimes_F (F^{(dk)}
\otimes_{F^{(k_+)}} L)$. Thus
\[\GM_1 = \{\gm\in\GM | dk/d'\text{ divides }
\val(\det(\gm))\}.\]
 The quotient $\GM_1/\GM'$ is
cyclic of order $d'$, generated by the projection of any element
$\gm_0\in\GM_1$ such that $\val(\det(\gm_0))=dk/d'$. Fix such a $\gm_0$. Then
\begin{equation}
\label{points}
 X(L) = \{P'\in X'(L^{(d')}) | \gm_0(P') =P'\}.
\end{equation}
These considerations prove the following
\begin{lemma}
\label{delpoints}
 $X(L)\neq\emptyset$ if and only if there exists
$P\in\oom{d}{F}(\overline{L})$ such that for every
$\sigma\in\Gal(\overline{L}/L)$ there exists an
element $\phi(\sigma)\in\GM_1$ satisfying
$\sigma(P)=\phi(\sigma)(P)$ and
$\phi(\sigma)\in\gm_0^i\GM'$ if
$\sigma_{|L^{(d')}}=Fr^i_L$. (Here $\GM_1$ acts on $\oom{d}{F}(\overline{L})$
through its projection to $\PGL_d(F)$.)
\end{lemma}

We return to the proof of the Theorem:

a) Assume that $X(L)\neq\emptyset$, and let $P$ be
as in the lemma. Set $\Dt:=\{\gm\in\GM_1|
\gm(P)=\tau(P)\text{ for some }
\tau=\tau(\gm)\in\Gal(\overline{L}/L)\}$. Then $\Dt$ is a subgroup of $\GM_1$,
containing $Z\GM$,
so it will suffice to show that $\Dt$ satisfies the conditions (i)-(iii) of the
Theorem.

(i) Since the natural $\GL_d(F)$-equivariant map
from $\oom{d}{F}$ to the Bruhat-Tits building
$\C{B}_F^d$ of $\PGL_d(F)$ is constant on the
Galois orbits,  $\Dt$ stabilizes a certain point of
$\C{B}_F^d$. By assumption B) on $\GM$, the group
$P\Dt \subset \PGL_d(F)$ is compact and discrete.
Hence it is finite.

(ii) By Lemma~\ref{delpoints}, the existence of $P$
forces $\Dt$ to contain an element from $\gm_0\GM'$.
Since $\Dt\supset Z\GM$, the statement follows.

(iii) Let $\dt_0$ be an element from
$\gm_0\GM'\cap\Dt\subset\GM_1$. Multiplying $\dt_0$
by an element of $Z\GM$ we may and will assume that
$\val_F(\det(\dt_0))=\val_F(\det(\gm_0))=kd/d'$.
Let $n$ be the order of the image of $\dt_0$ in
$P\Dt$, and let $\dt_0^n=a\in F\m \subset L\m$.
Then $\det(\dt_0)^n= a^d$, hence $\val_L(a)=n
\val_L(\det(\dt_0))/d=ekn/d'$.

Let $L'$ be the field of rationality of $P$
(corresponding by Galois theory to the stabilizer
$H\subset\Gal(\overline{L}/L)$
of $P$). Then for every $\sigma\in\Gal(\overline{L}/L)$ and every $\tau\in H$ we have
$$\tau(\sigma(P))=\tau(\phi(\sigma)(P))=\phi(\sigma)(\tau(P))=\phi(\sigma)(P)=\sigma(P),$$
hence $L'$ is a Galois extension of $L$. Also we have a natural surjective homomorphism,
$\pi:\Dt\to\Gal(L'/L)$ such that $\pi(\dt)(P)=\dt(P)$ for each $\dt\in\Dt$.

Let $v\in (L')^d$ be a representative of
$P\in\oom{d}{L}(L')\subset\B{P}^{d-1}(L')$. Then by
our assumption, there exist $\tau\in \Gal(\overline{L}/L)$ and $\lambda\in L'$ such that
$\dt_0(v)=\lambda \tau(v)$. Since the action of $\dt_0\in\GL_d(F)$ commutes with that of $\tau$,
we get $av=\dt_0^n(v)=\lambda
\tau(\lambda) \tau^2(\lambda)
...\tau^{n-1}(\lambda)(v)$. Hence $a=\lambda
\tau(\lambda) \tau^2(\lambda)
...\tau^{n-1}(\lambda)$. Taking valuations, we get
$\val_L(\lambda)=\val_L(a)/n=ek/d'$. Since
$\val_{L'}(\lambda)$ is an integer, $d'|ek
e(L'/L)$, so it will suffice to show that $e(L'/L)$
divides $m$.

But $e(L'/L)$ obviously equals the ramification
degree of $L'L^{(d')}$ over $L^{(d')}$, hence it
divides  the degree $[L'L^{(d')}:L^{(d')}]$. This
degree is equal to the number of conjugates of $P$
over $L^{(d')}$. But conjugates of $P$
over $L^{(d')}$ form a homogeneous space for the action of the group $\Dt'/Z\GM$, hence the number
of conjugates divides the order of  $\Dt'/Z\GM$, which is $m$.
This completes the proof of part a).

b) We will show the existence of local points by a
case-by-case analysis.

I) Assume first that $d|efk_+$ or, equivalently, that
$d'|ek$. Let $\dt_0$, $n$ and $a$ be as in the
proof of iii) in a). It will suffice us to show
that there exists a point $P\in\oom{d}{F}(L^{(n)})$
such that $\dt_0(P)=Fr_L(P)$.

The assumption $d'|ek$ implies that $n|\val_L(a)$.
Therefore there exist $\lambda\in L^{(n)}$ such
that $N_{L^{(n)}/L}(\lambda)=a$. Hence
$\dt':=\lambda^{-1}\dt_0\in\GL_d(L^{(n)})$
satisfies
$\Fr_L^{n-1}(\dt')\Fr_L^{n-2}(\dt')...\dt'=1$. By
Hilbert Theorem 90 for $\GL_d$, there exists
$B\in\GL_d(L^{(n)})$ such that $\dt'=\Fr_L(B)\cdot
B^{-1}$. Hence for every vector $v\in L^d$, the
image $P=P_v$ of $Bv\in(L^{(n)})^d$ in
$\B{P}^{d-1}(L^{(n)})$ satisfies
$\dt_0(P)=\Fr_L(P)$. Therefore it remains to show
that there exists  $v\in L^d$ such that the
corresponding $P_v$ belongs to $\oom{d}{F}$.

Let $W$ be the $L^{(n)}$\/-vector subspace of
$\Mat_d(L^{(n)})$ spanned by the Galois conjugates
of $B$ over $L$. Then $W$ is invariant under the
action of the Galois group $\Gal(L^{(n)}/L)$, hence
is defined over $L$ by descent theory. Therefore
$W$ contains an element $B'\in\GL_d(L)$.

By our assumption, $d|efk_+$, hence
$[L:F]=efk_+\geq d$. In particular, $L$ contains
$d$ elements which are linearly independent over
$F$. Equivalently, there exists $v'\in L^d$ not
lying in any $F$-rational hyperplane. We claim that
for $v:=(B')^{-1}(v')$ the corresponding $P_v$ lies
in $\oom{d}{F}$. In fact, let  $w\in F^n$ be a  row
vector such that $wBv=0$. Since $w$ and $v$ are
$L$-rational, we get $wCv=0$ for every matrix $C\in
W$. In particular, we have $wv'=wB'v=0$. By our
choice of $v'$, the vector $w$ is therefore the
trivial vector, so that $P_v\in\oom{d}{F}$.

II) Assume now that $d$ is prime but it does not divide
$fek_+$. Then $d'=d$, $k_+=k$ and $[\Dt:\Dt']=d$.
Let $\dt_0$, $ n$ and $a$ be as above.

We claim that for every field extension $\wt{L}/F$
whose ramification degree $e(\wt{L}/F)$ is prime to
$d$ (in particular for $\wt{L}=L$), the subalgebra
$\wt{L}_0=\wt{L}[\dt_0]\subset\Mat_d(\wt{L})$ is a
totally ramified field extension of $\wt{L}$ of
degree $d$. Indeed, since the minimal polynomial of
$\dt_0$ divides $x^n-a$, the algebra $\wt{L}_0$ is
either a field or a direct sum of fields. Let
$\wt{L}'_0$ be one of the direct factors of
$\wt{L}_0$, and let $\dt'_0$ be the image of
$\dt_0$ in $\wt{L}'_0$. Then $(\dt'_0)^n =a$,
therefore
$\val_{\wt{L}'_0}(\dt'_0)=\val_{\wt{L}'_0}(a)/
n=e(\wt{L}'_0/F)ekn/nd= e(\wt{L}'_0/F)ek/d$.
Since $\val_{\wt{L}'_0}(\dt'_0)$ is an integer,
our assumption implies that
$d|e(\wt{L}'_0/\wt{L})|[\wt{L}'_0:\wt{L}]\leq
[\wt{L}_0:\wt{L}]=d$. From this the statement
follows.

It follows from our assumption that $d^2$ divides
the order of $P\Dt$. We distinguish two cases: i)
$d^2|n$ and ii) $d^2$ does not divide $n$.

Case i) We will prove that there is a cyclic
extension $M/L$ of degree $n$, whose ramification
degree is $d$ such that $a\in \Nm_{M/L}M\m$.
Write $n$ as a product $d^t n'$, where $n'$ is
prime to $d$. As before,
$\val_L(a)=ekn/d=d^{t-1}ekn'$, hence
$n'|\val_L(a)$. It follows that $a$ belongs to
$\Nm_{L^{(n')}/L}(L^{(n')})\m$. If we find a cyclic
extension $M'/L$ of degree $d^t$ with ramification
degree $d$ such that  $a\in \Nm_{M'/L}(M')\m$,
then the composite field $M:=M'L^{(n')}$ satisfies
the required property. Indeed, as $M'$ and $L^{(n')}$ are linearly disjoint over $L$,
 we have
$\Nm_{M'L^{(n')}/L}(M'L^{(n')})\m=\Nm_{M'/L}(M')\m
\cap\Nm_{L^{(n')}/L}(L^{(n')})\m$.

By Local Class Field Theory
to construct $M'$  it is equivalent to construct an
open subgroup $H\subset L\m$, containing $a$ and
contained in $\{l\in L\m| d^{t-1} \text{ divides }
\val_L(l)\}$ such that $L\m/H\cong\B{Z}/d^t\B{Z}$.

First we show that $a$ is not a $d^{\text{th}}$ power in
$L$. In fact, suppose that $a=b^d$ for some $b\in L$.
Let $\eta$ be a primitive $d^{\text{th}}$ root of unity
inside $\overline{L}$. Then the ramification degree
of $\wt{L}:=L(\eta)$ over $L$ (hence over $F$) is prime to $d$, so
by the claim above the algebra
$\wt{L}_0=\wt{L}[\dt_0]\subset\Mat_d(\wt{L})$ is a
field. The equality
$(\dt_0^{n/d}-b)(\dt_0^{n/d}-\eta
b)...(\dt_0^{n/d}-\eta^{d-1}b)= \dt_0^n-a=0$ then
implies that some $\dt_0^{n/d}-\eta^i b$ equal to
$0$. Hence $\dt_0^{n/d}$ is central in
$\Mat_d(\wt{L})$, contradicting the fact that $\dt_0$ has an order $n$ modulo center.

The group $A:=L\m/(L\m)^{d^t}$ is a finite abelian $d$-group. Let $\bar{a}$
and $\bar{\pi}$ be the images in $A$ of $a$ and of
some uniformizer of $L$ respectively, and set
$A':=\C{O}_{L}\m/(\C{O}_L\m)^{d^t}\subset A$. We claim
 that there exists
a subgroup $H'\subset A'$ such that $A$ decomposes
as a direct sum
$\langle\bar{\pi}\rangle\oplus\langle\bar{a}\rangle
\oplus H'$, where by $\langle\bar{\pi}\rangle$
(resp. $\langle\bar{a}\rangle$) we denote the
cyclic subgroup generated by $\bar{\pi}$ (resp.
$\bar{a}$). Indeed, using
the fact that $d^{t-1}|\val_L(a)$ and that $a$ is
not a $d^{\text{th}}$ power in $L$, we first see that
that cyclic subgroups $\langle\bar{\pi}\rangle$ and $\langle\bar{a}\rangle$
have a trivial intersection.
Put $a':=a\pi^{-\val_L(a)}\in\C{O}_L\m$. It remains to show that
$\langle\bar{a'}\rangle$ is a direct summand in $A'$, or equivalently
that $a'$ is not a  $d^{\text{th}}$ power. As $d|\val_L(a)$ and $a$ is
not a $d^{\text{th}}$ power in $L$, the statement follows.
Now we can take  $H\subset L\m$ be the inverse image of
$\langle\bar{a}\rangle\oplus H'$ in $L\m$.

Let $\sigma$ be a generator of $\Gal(M/L)$ such
that $\sigma_{|L^{(n/d)}}=\Fr_L$. As in the
previous case, it will suffice to find a point
$P\in\oom{d}{F}(M)$ such that $\dt_0(P)=\sigma(P)$.
Let $\lambda$ be an element of $M$ such that
$\Nm_{M/L}(\lambda)=a$. As before, there exists
$B\in\GL_d(M)$ such that
$\sigma(B)=\lambda^{-1}\dt_0 B$. Then for every
$v\in L^d$ the image $P=P_v\in\B{P}^{d-1}(M)$ of
$Bv\in M^d$ satisfies $\dt_0(P)=\sigma(P)$, so it
remains to show that some $Bv$ does not lie in an
$F$-rational hyperplane. In fact, suppose that
$wBv = 0$ for some row vector $w\in F^d$.
Since $\sigma(Bv)=\lambda^{-1}\dt_0 Bv$,
it follows that $w\dt_0^i Bv = 0$ for all $i$.
Since the minimal polynomial of $\dt_0$ has degree
$d$, we get that for a generic $v$ the vectors $\dt_0^i Bv=B(B^{-1}\dt_0 B)^i v$
span all of $M^d$. For such a $v$ we therefore get $w=0$, as claimed.

Case ii) Let $\bar{\dt}_0$ be the image of
$\dt_0^{n/d}$ in $P\Dt$, and let $\Dt_d$ be a
$d$-Sylow subgroup of $P\Dt$ containing
$\bar{\dt}_0$. By our assumptions, $\Dt_d\cap P\Dt'$ is a non-trivial $d$-subgroup,
hence it contains a central element $\bar{\dt}_1$ of $\Dt_d$ of order $d$. 
Since $n/d$ is prime to $d$,
we get that $\bar{\dt}_0\notin P\Dt'$. Therefore $\bar{\dt}_0$ and $\bar{\dt}_1$ generate a
subgroup, isomorphic to $(\B{Z}/d\B{Z})^2$.

Let $\dt_1\in\Dt$ be a representative of
$\bar{\dt}_1$, then $(\dt_1)^d\in F\m$. Replacing $\dt_0$ by $\dt_0^{n/d}$, we get
$\dt_0\dt_1=b\dt_1\dt_0$ for some $b\in F\m$.
Taking the determinants we see that $b^d=1$. We
claim that $b\neq 1$. In fact, if $b$ would equal
to $1$, then $\dt_1$ would belong to the
centralizer of $\dt_0$ in $\Mat_d(F)$, which is
$F(\dt_0)$. Since $F(\dt_0)$ is a totally ramified
extension of $F$, the only elements of $F(\dt_0)\m$
whose $d^{\text{th}}$ power is in $F\m$ are of the form
$c(\dt_0)^i$ for some $c\in F\m$ and some $i$. But
then we get $\bar{\dt}_1=(\bar{\dt}_0)^i$,
contradicting to our assumptions. Hence
$b$ is a primitive $d^{\text{th}}$ root of unity.

The characteristic polynomial of $\dt_0$ is
irreducible over $F$, therefore $\dt_0$ has $d$ distinct
eigenvectors $V_1,...V_d$ (conjugate over $F$).
They correspond to $d$ conjugate fixed points
$P_1,...,P_d\in\oom{d}{F}$ of $\dt_0$. To finish
the proof of the theorem it will suffice to show
that the cyclic group generated by $\bar{\dt}_1$
acts simply-transitively on the $P_i$'s.

Let $\lambda_1,...\lambda_d$ be the eigenvalues of
$\dt_0$, corresponding the $V_i$'s, that is, $\dt_0
(V_i)=\lambda_i V_i$ for all $i=1,...,d$. Then for
each $i=1,...,d$ and each $j=1,...,d-1$ we have
$\dt_0(\dt_1)^j
(V_i)=b^j(\dt_1)^j\dt_0(V_i)=b^j\lambda_i(\dt_1)^j
(V_i)$.  So that $(\bar{\dt}_1)^j (P_i)\neq P_i$,
as was claimed.
\end{pf}

\begin{remarks}
{\rm a) By similar analysis, we can show that converse
of a) holds in some more cases, but we don't know
whether it holds in general.

b) We suspect that there should be a simpler
``topological'' proof of Theorem~\ref{T:Main},
which uses the fact that the Bruhat-Tits building
of $\PGL_d(F)$ can be naturally embedded into
(Berkovich's) $\oom{d}{F}$.

c) The curve case of Theorem~\ref{T:Main} was
studied in \cite{JL1}, and generalizing these
results to the higher dimensional case is one of
our aims here --- see also the next section. It is
in fact possible to obtain our results here through
a generalization of the method of loc.\ cit.,
namely an analysis of the special fiber and the
resolution of its singularities.}
\end{remarks}

\begin{Cor} \label{C:torsion}
If $P\GM$ is torsion-free, then $X(L)\neq\emptyset$
if and only if $dk|fk_+$.
\end{Cor}

\begin{pf}
Since  $\val_F(\det(Z\GM))=dk\B{Z}$, the corollary
follows immediately from the theorem.
\end{pf}

\section{Local Points on Shimura Curves}
We shall now specialize to the case of curves,
which is Case~1 of Section~\ref{padic} with $r=1$.
We restrict the level at the prime $\cP$ of $F$
corresponding to $v_1$ by taking a compact open
subgroup $S$ of $\cG = \cG'$ and considering a
geometrically connected component $X$ of
$S\backslash\Shimp$. By Theorem~\ref{fieldofdef},
$X$ is defined over $F_\cP^{(k_+)}$. Put
$\oom{2,\nr}{F_\cP} = \oom{2}{F_\cP}\hat{\otimes}_{F_\cP}
F_\cP^\nr$. By Theorem~\ref{int},  $X$ is
$F_\cP^{(k_+)}$\/-isomorphic to
$\Gamma(gSg^{-1})\bs \oom{2,\nr}{F_\cP}$ for some $g\in \cG$, where we write $\Gamma(T)$
instead of $T\cap G(F)$. We also
restrict the level away from $\cP$ by assuming that
\begin{itemize}
\item[(*)] the norm of $S$ and its
intersection with the center agree as
subgroups of $(\bbA_\rf^{\cP})^\times$.
\end{itemize}

Our restriction of the level at $\cP$ means that
the level (of $G^\rint$\/) at $\cP$ is the maximal
compact subgroup compact modulo the center. (Some
cases of larger level subgroups at $\cP$ are
discussed in \cite{JL3}, and a little bit is known
about the simplest case of smaller level at $\cP$
(\cite{Tei}), but further information seems to be
required to handle the general case.) Assumption
(*), on the other hand, is mainly for convenience.

In  the notation of Section~\ref{concomp}, $T = \Nm
S (= Z_\rf^{\cP} \cap S)$, and ${\wt{T}} =
T\cO_{F_\cP}^\times$. We have $\bM =
\Res_{F/\bbQ}\Gm$, so that $\Cl_{\wt{T}}^+ = {\wt{T}}
F^+_\infty
     \backslash \bbA^\times/ F^\times$.
Here $F^+_\infty$ is the subgroup of
everywhere positive elements in
$F_\infty^\times =
      \prod_i F_{\infty_i}^\times$.
Define also
$\Cl_{\wt{T}} = {\wt{T}} F^\times_\infty
     \backslash \bbA^\times /F^\times$,
and let $k$ be the order of the image
$\varphi$ of $\pi_\cP$ in $\Cl_{\wt{T}}$. We
view the center $Z\GM(S)$ and
$\Nm\Gamma(S)$ as subgroups of $\cOstFP$,
the units away from $\cP$ of $F$. Then we have
the following
\begin{lemma}
\label{kkplus}
a) \ $\Nm\Gamma(S) = T \cap F^\times_+$, with
$F^\times_+$ be the totally positive elements of
$F$.\\[0.05in]
b) \ $k\bbZ = \val_\cP Z\GM(S)$ and
 $k_+\bbZ = \val_\cP \Nm\Gamma(S)$.
\end{lemma}
\begin{pf}
$\Nm\Gamma(S) \subseteq T \cap F^\times_+$
holds since $B$ is definite. Conversely, take
$t \in T \cap F^\times_+$ (viewed in
$(\bbA^{\cP}_{\rf})\m $). By the Hasse
principle $t = \Nm b$ for some
$b\in B^\times$, and also $t = \Nm s$ for some
$s\in S$. Then $\Nm bs^{-1} = 1$ in
$G_\rf^{\cP}$. By the Eichler-Kneser strong
approximation theorem, $bs^{-1} = b_1s_1$ for
some $b_1 \in B$ and $s_1 \in S$, both of
norm $1$. Then $b_1^{-1}b =  s_1s$ is in
$\Gamma(S)$ and has norm $t$, proving Part
a).
Part b) is equivalent to the exactness of
the two sequences
\[Z\GM(S) \sra{\alpha}
  F_\cP^\times/\cO_{F_\cP}^\times \sra{\beta}
  \Cl_{\wt{T}} \qquad \mbox{ and } \qquad
  \Gamma(S) \sra{\alpha'}
  F_\cP^\times/\cO_{F_\cP}^\times \sra{\beta'}
  \Cl_{\wt{T}}^+ \,.\]
Here $\alpha$ takes the $\cP$\/-component,
$\alpha'$ takes the $\cP$\/-component of the norm,
and $\beta$, $\beta'$ are the natural maps. The
exactness is routine, and left to the reader.
\end{pf}

We can now apply the results of
Section~\ref{locpoints}. As $d=2$, we have two possibilities: either
$k_+ = k$, in which case
$X$ is the quadratic unramified twist of $X'$ by
$w_p(=$ the image of $\gm_0$ in $\Aut X'$) or
$k_+ = 2k$ in which case $X\simeq X'$. In the latter case we will say that $X$
is  {\em Mumford-uniformized}\/, even though
strictly speaking $X$, $X'$ are obtained from a
variety which is Mumford-uniformized over $F$ by
extending scalars to $F^{(k_+)}$\/.

For a finite set $W$ of finite primes
$\neq \cP$ of $F$, we denote by $K(W)$ the
principal congruence subgroup of (squarefree)
level $\prod_{\cQ \in W} \cQ$ in
$G_\rf^{\cP}$. We will now show that for a
``sufficiently small'' $S$ we have $k_+ = k$:
\begin{proposition}
\label{small}
Let $W_1$ be a finite set of primes of $F$.
Then there exists a finite set $W_2$ of
finite primes of $F$, disjoint from
$W_1 \cup \set{\cP}$, such that $k_+ = k$ for
any open subgroup $S \subset K(W_2)$.
\end{proposition}
\begin{pf}
By the Unit Theorem, we have
$\cOstFP \simeq (\bbZ/2\bbZ) \times \bbZ^g$.
For $0 \leq i \leq g$, choose
$\epsilon_i \in \cOstFP$ whose classes modulo
squares are a basis for the $\bbF_2$\/-vector
space $\cOstFP/(\cOstFP)^2$. A routine
application of Chebotarev's Theorem shows
that there are primes $\cQ_i$,
$0 \leq i \leq g$, not in
$W_1 \cup \set{\cP}$, such that $\epsilon_i$
has a non-square residue modulo $\cQ_j$ if and
only if $i = j$. Put $W_2 = \set{\cQ_i}$.
Then for any open subgroup $S \subset K(W_2)$
we have
$Z\GM(S)=Z(S) \cap F^\times \subset (\cOstFP)^2$,
so a fortiori
$Z\GM(S)= \Nm\GM(S)$.
The Proposition follows.
\end{pf}
\begin{remarks} {\rm 1.\ It follows from the proof
that for $S \subset K(W_2)$ as above, the
$\cP$\/-adic valuation of each
$x \in Z\GM(S)$ is even, hence
$k_+ = k$ is even.\\[0.02in]
2.\ The Proposition can be strengthened in several
ways. If $\epsilon_1, \dots, \epsilon_h$ are in the
subgroup $\cO_{F,\cP}^+$ of the totally positive
elements of $\cOstFP$ and generate it modulo
squares, it suffices to take $W_2 = \set{\cQ_i}$,
$h< i \leq g$. It is also possible to choose $W_2$
in a set of positive density (with the obvious
definition of positive density here), not merely to
avoid $W_1 \cup \set{\cP}$. In addition, if the
primes of $W_1 \cup \set{\cP}$ are all prime to
$2$, one can replace $K(W_2)$ by an open subgroup
of $G_\rf^{ \cP}$ whose level divides some power of
$2$. We shall not need these facts in the
sequel.\\[0.02in]
3.\ The Proposition and its proof clearly
generalize from $d=2$ to any $d$.}
\end{remarks}

We will need the following
\begin{lemma} \label{L:tree}
\label{torlem} Let $\gm\in \GL_2(F_\cP)$ has
finite order modulo the center and suppose that
$\val_\cP \det\gm$ is odd. Then $\gm$ reverses a
unique edge $e_\gm$ of the Bruhat-Tits tree $\C{B}_F^2$.
\end{lemma}
\begin{pf}
The subgroup of $\PGL_2(F_\cP)$ generated by $\gm$
is contained in a maximal compact subgroup, hence it
fixes a vertex or an edge. Since $\gm$ moves each
vertex to an odd distance, it cannot fix any vertex,
hence must reverse a unique edge.
\end{pf}
Theorem~\ref{T:Main} now says the following:
\begin{theorem}
\label{main} a) If $k_+ = 2k$, or if $f$ is even,
then $X(L) \neq \emptyset$.

b) Suppose $k_+ = k$
and $f$ odd. Then $X(L) \neq \emptyset$ if and only
if there exists an element $\gm \in \GM(S)$, of finite
order modulo the center satisfying $\val_\cP \det
\gm = k$, such that either $ek$ is even, or $ek$ is
odd and the stabilizer $H$ of $e_\gm$ in
$P\Gamma(S)$ as an {\em oriented} edge has even
order.
\end{theorem}
\begin{pf}
Since $d=2$ is a prime, part a) follows immediately from
Theorem~\ref{T:Main} by taking $\Dt=Z\GM$. To deduce part b), suppose
first $X(L) \neq \emptyset$. Then the subgroup
$P\Dt$ of Theorem~\ref{T:Main} contains an element
$\gm$ as needed, and either $ek$ is even or
$P\Dt'$ has even order. But if $ek$ is odd, then
$P\Dt'$ is precisely the stabilizer of $e_\gm$ as
an oriented edge. Conversely, if $ek$ is even let
$\Dt$ be the subgroup of $P\GM$ generated by $\gm$;
and if $ek$ is odd let $\Dt$ be the stabilizer in
$\GM$ of the unoriented edge $e_\gm$, defined in Lemma~\ref{L:tree}. Then it is
immediate that $\Dt$ has the properties required in
Theorem~\ref{T:Main}, so $X(L) \neq \emptyset$.
\end{pf}

If $S$ is small enough, then $P\Gamma(S)$ is
torsion free. If $S$ moreover satisfies the
condition of Proposition~\ref{small}, we shall
say we are in the {\em asymptotic  case}. It
then follows from the Theorem that
\begin{equation*}
\mbox{$k_+ = k$, and $X(L) \neq \emptyset$ if
and only if $f$ is even.}\tag{ASYMP}
\end{equation*}

If we only assume that $P\Gamma'(S)$ is
torsion free, the Theorem has the following
\begin{corollary}
If $P\Gamma'(S)$ is torsion free,
then $X(L) = \emptyset$ if and only if
$k_+ = k$, $f$ is odd, and either
$ek$ is also odd or $P\Gamma(S)$ is torsion
free.
\end{corollary}
\begin{pf}
By the theorem, we may assume that $k_+ = k$, that $f$ is odd
and that $P\Gamma(S)$ has a torsion. It remains to show that
 $X(L)=\emptyset$ if and only if $ek$ is odd. Assume first that
$ek$ is odd.  As  $P\Gamma'(S)$ is torsion free, no non-trivial element of $P\GM(S)$
stabilize an oriented edge, so we get  $X(L)=\emptyset$. Assume now that
$ek$ is even. Let $\gm'\in\GM$ be any torsion moduli
center element. Since $P\Gamma'(S)$ is torsion free,  $\val_\cP \det
\gm'\equiv k\pmod{2k}$. Therefore modifying $\gm$ by an element
from $Z\GM(S)$, we get an element $\gm$ with $\val_\cP\det\gm=k$, as claimed.
\end{pf}

At the other extreme, we can pin down the
situation rather precisely when
$X(L) \neq \emptyset$  and $fek_+$ is odd (so
that $k = k_+$\/). This will require two lemmas,
which we state in slightly more general form
than necessary.

For an $n^{\text{th}}$ root $\zeta$ of $1$, let
$\bbQ(\zeta)^+$ be the maximal totally real
subfield of the field $\bbQ(\zeta)$ of
$n^{\text{th}}$ roots of $1$. Then we have the
following

\begin{lemma}
\label{valp}
Let $\zeta$ be a primitive $n^{\text{th}}$ root of
$1$, with $n>2$. Let $\cQ^+$ be a prime ideal
of $\bbQ(\zeta)^+$ of residue characteristic
$q$. Then
$\val_{\cQ^+} (1+\zeta)(1+\ov{\zeta}) = 1$
if $n = 2q^r$ for  $r\geq 1$. Otherwise
this valuation is $0$.
\end{lemma}
\begin{pf}
Let $A$ be the positive integer
$\Nm_{\bbQ(\zeta)/\bbQ}(1+\zeta)$. We claim
that $A = \ell$ if $n = 2\ell^r$ for a prime
$\ell$, and $A = 1$ otherwise. Indeed, write
$n = 2^rt$ with $t$ odd and $r \geq 0$. If
$r = 0$ then $A = 1$ because
$\Nm_{\bbQ(\zeta)/\bbQ}(1-\zeta) =
       \Nm_{\bbQ(\zeta)/\bbQ}(1-\zeta^2)$.
For any $u$, let $\Phi_u$ denote the $u^{\text{th}}$
cyclotomic polynomial. If $r = 1$ then
$A = \prod_\eta (1-\eta)$,
where $\eta$ goes over all primitive $t^{\text{th}}$
roots of $1$, so $A = \Phi_t(1)$. Then if $t$
is a power of a prime $\ell$, we have
$\Phi_t(x) = (x^t - 1)/(x^{t/\ell} - 1) =
x^{(\ell - 1)t/\ell} + \dots + 1$, so $A =\ell$;
else $t = t't''$ with $t'$, $t''$ odd and
relatively prime. Then $\Phi_t(x)$ divides
$Q(x) = (x^{(t' - 1)t''} + \dots + 1)
     /(x^{t' - 1} + \dots + 1)$,
so $A|Q(1) = 1$. Similarly, if $r \geq 2$,
then $A = \Phi_n(1)$. Then if $t = 1$, we have
$\Phi_n(x) = (x^{n/2} + 1)$, so $A = 2$,
whereas if $t > 1$, then $\Phi_n(x)$ divides
$Q(x) = \Phi_{2^r}(x^t)/\Phi_{2^r}(x)
      = (x^{n/2} + 1)/(x^{2^{r-1}} + 1)$,
so $A|Q(1) = 1$. This proves our claim in all
cases. The Lemma is now immediate if
$n\neq 2q^r$, $r > 0$. Else, recall that
$1 + \zeta$ is a uniformizer for the (totally
ramified) prime $\cQ$ of $\bbQ(\zeta)$ above
$q$ and $\cQ^+$. Then
$\val_{\cQ^+} (1+\zeta)(1+\ov{\zeta})
   = \val_{\cQ} (1+\zeta) = 1$,
proving the Lemma.
\end{pf}

Assume that $k_+=k$, and let $\gamma$ be as in Theorem~\ref{main}. Then
$F(\gamma)$ is a CM extension of $F$, since it
splits the definite algebra $B$ and is $\neq F$.
Then $\Nm \gamma$ is a (totally positive) generator
of $\cP^{k}$. Writing $\Nm$ for $\Nm_{K/F}$ for a
quadratic extension $K$ of $F$ is unambiguous,
because $\Nm_{K/F}$ and the reduced norm of $B$
agree for any $F$\/-embedding of $K$ into $B$.
Likewise we let $\ov{x}$ denote both the main
involution of $B$ and the (complex) conjugate of
any element $x\in K \subset B$.
\begin{lemma}
\label{rootof1}
Let $F(\gamma)$ be a {\rm (}quadratic\/{\rm )}
CM extension of $F$, such that some power
$\gamma^m$ is in $F$, with $m$
{\rm (}$>1${\rm )} minimal. Suppose that
$\Nm \gamma$ generates a
power $\cP^k$, $k>0$, of a prime ideal
$\cP$ of $F$, of residue characteristic
$p$. Set $\zeta = \gamma/\ov{\gamma}$.
Then $\gamma$ is an integer of $F(\gamma)$,
$\zeta$ is a primitive $m^{\text{th}}$ root of
$1$, and if $m$ is even, then
$\gamma^m$ is totally negative. Moreover,
\begin{enumerate}
\item[(i)] If $m=2$ then $\gamma^2 = -u$, where
$u\in F$ is a positive generator of $\cP^k$.
\item[(ii)] If $m\neq 2$ then $F(\gamma) = F(\zeta)$, and
$\gamma = s(1+\zeta)$ for some $s \in F^\times$.
\end{enumerate}
\end{lemma}
\begin{pf} $\gamma$ is an algebraic integer
since its power is. The order of $\zeta$ is the the smallest positive integet $n$ such that
$\gamma^n = \ov{\gamma}{}^n$ or equivalently that $\gm^n\in F$. Thus $n=m$, as claimed.
If $m$ is
even, then $\zeta^{m/2} = -1$. Hence
$\gamma^m = -\gamma^{m/2}
   \ov{\gamma}{}^{-m/2} = -\Nm \gamma^{m/2}$
is totally negative as asserted, and we also
get (i). For (ii), define
$s = \gamma\ov{\gamma}/
                      (\gamma + \ov{\gamma})$,
which makes sense since  $\ov{\gamma} \neq -\gamma$.
Then $\gamma = s(1+\zeta)$ as asserted.
\end{pf}
\begin{remark}
{\rm $s$ need not belong to $\cOstFP$: take
$F = \bbQ(\sqrt{6}\,)$, $\cP = (3,\sqrt{6}\,)$,
and $\gamma = \sqrt{6}(1+\sqrt{-1}\,)/2$. Then
$m =4$, and $\gamma = \sqrt{3}\zeta_8$, but
$\sqrt{6}/2$ is not in $\cO^\times_{F,\cP}$.}
\end{remark}

 Let $t(F)$ be the maximal
integer such that the maximal totally real
subfield $\bbQ(\zeta_{2^{t(F)}})^+$ of
$\bbQ(\zeta_{2^{t(F)}})$ is a subfield of $F$. Then $t(F)\geq 2$.
Let $\tau$ be the unique prime of $\bbQ(\zeta_{2^{t(F)}})^+$
lying (and totally ramified) above $2$.
We now have the following
\begin{proposition}
\label{crit}
In the notation of Theorem~\ref{main}, suppose
that $fek_+$ is odd. Then $X(L)\neq\emptyset$ if and only if either
of Conditions (a) or (b) below holds:
\begin{enumerate}
\item[(a)] $\cP$ lies above $\tau$ with odd
ramification index $e(\cP/\tau)$, and there
exists an element $\gamma'$ of $\Gamma(S)$,
of order $2^{t(F)}$ modulo the center, such
that $\val_\cP \Nm(\gamma') = k$.
\item[(b)] There exist elements $\alpha$,
$\gamma'$ in $\Gamma(S)$, of order $2$ modulo
the center, such that $\val_\cP \Nm(\gamma') = k$,
$\val_\cP \Nm(\alpha) = 0$, and
$\alpha\gamma = -\gamma\alpha$.
\end{enumerate}
\end{proposition}
\begin{pf}
Suppose first that (a) or (b) hold. Lemma~\ref{kkplus}
and the fact that $f$ is odd enable us to modify
$\gamma'$ by a central element in $Z\GM(S)$ and
hence to replace it by and element $\gamma$
with the same order modulo the center and such that
$\val_\cP \Nm(\gamma) = fk_+$. Since $fk_+$
is odd, $\gamma$ reverses an edge.
The stabilizer $H$ of this edge contains
$\ov{\gamma}\vphantom{\gamma}^2$ in case~(a),
and  $\ov{\alpha}$ in case~(b). Hence
its order $|H|$ is even in either case. By
Theorem~\ref{main}, $X(L) \neq \emptyset$.

Conversely, if $X(L) \neq \emptyset$ let $\gamma$
be as in Theorem~\ref{main}. Then an odd power of
$\gamma$ is of order $2^t$, $t \geq 1$, modulo the
center. Modifying this power by a central element
of $\Gamma(S)$ we get an element $\gamma_1$, of
order $2^t$ modulo the center, such that
$\val_\cP\Nm\gamma_1 = fk_+$. We now replace $\gamma$
by $\gamma_1$. The new $\gamma$ has the same fixed
edge $d'$ as before, and satisfies the properties
of Theorem~\ref{main}, and its order modulo the
center is $m=2^t$.

Suppose first $t=1$.
Theorem~\ref{main} gives an element $\beta\in H$,
of order $2$ modulo the center.
Let $A$ be the subgroup of $P\Gamma(S)$ generated
by the images $\ov{\beta}$ and $\ov{\gamma}$ of $\beta$ and $\gamma$.
Since $A$ acts discretely on $\C{B}_F^2$ fixing the unoriented edge
$\set{d',\ov{d'}}$, it is a finite group. Since $A$
is generated by $2$ involutions it is a dihedral
group, and the cyclic subgroup $A_0$ generated by
$\ov{\beta}\ov{\gamma}$ is of index $2$. Also $\ell = |A_0|$
is even, because $A$ has {\em another} subgroup
$A_1 \neq A_0$ of index $2$, namely the elements of
$A$ fixing $d'$. Hence the group generated by
$\ov{\gamma}$ and $\ov{\alpha}$, where $\alpha =
(\beta\gamma)^{\ell/2}$, is a non-cyclic group of
order $4$. Replacing $\alpha$ by $\alpha\gamma$ if
necessary, we may assume that $\alpha d' = d'$. In
$\Gamma(S)\subset B^\times$ this implies that
$\alpha^2$ and $\gamma^2$ are in $F^\times$, that
$\alpha$ is in $\Gamma'(S)$, and that
$\gamma\alpha = u\alpha\gamma$ for some $u\in
F^\times$. Taking norms (to $F$) we see that
$\Nm_{B/F}u = u^2 = 1$, so that $u = \pm 1$.
However $u = 1$ is impossible, since in this case $\alpha$ and
$\gamma$ would generate a commutative subgroup of
$\GL_2(\bbC)$ whose image modulo the center is
finite and non-cyclic, which cannot happen. This
gives case (b) of the Proposition.

Suppose next $t \geq 2$. By Lemma~\ref{rootof1},
$\gamma = s(1+\zeta_{2^t})$ and
$F(\gamma) = F(\zeta_{2^t})$, so that $F$
must contain $F_0 = \bbQ(\zeta_{2^t})^+$.
Then $\Nm\gamma = s^2 \theta$ where
$\theta = (1+\zeta_{2^t})(1+\ov{\zeta_{2^t}})$
generates the unique prime $\tau_0$ of $F_0$, lying over $2$. By Lemma~\ref{valp}, we get
\[ fk_+ = \val_\cP s^2 \theta \equiv
     e(\cP/\tau_0) \val_{\tau_0} \theta
        =  e(\cP/\tau_0) \pmod{2}. \]
This shows that $e(\cP/\tau_0)$ is odd.

Finally we check that $t = t(F)$. From
its definition $t \leq t(F)$; but if
$t(F) \geq t+1$ then $e(\cP/\tau_0)$ would
be divisible by the ramification index of
$\tau_0$ in $\bbQ(\zeta_{2^{t+1}})^+$, which
is $2$, contradicting the fact that
$e(\cP/\tau_0)$ is odd. Hence $t = t(F)$,
concluding the proof of the Proposition.
\end{pf}

\begin{remark}
\label{quaternion}
{\rm In case (b) above, the isomorphism type of
the quaternion algebras $B$ and $\Bint$
is almost forced. Indeed,  for
$a_1,a_2 \in F^\times$ let $B(a_1,a_2)$
denote the quaternion algebra over $F$
with basis $1$, $\qi$, $\qj$, $\qk$
satisfying $\qi^2 = a_1$, $\qj^2 = a_2$,
and $\qk = \qi\qj = -\qj\qi$. Then
$B(a_1,a_2)$ is ramified at a place $v$ of
$F$ if and only if the Hilbert symbol
$(a_1,a_2)_v$ is $-1$. In case (b) we
have $B \simeq B(-\Nm\gamma,-\Nm\alpha)$, so
that $B$, in addition to being totally
definite, is unramified away from primes
above $2$ (we know it is unramified at
$\cP$\/). For example, if $F = \bbQ$ we get
that $B \simeq \bbH = B(-1,-1)$, the
Hamilton quaternions over $\bbQ$; if $F$ is
(real) quadratic, then $B$ is isomorphic to
the base change of $\bbH$ to $F$, or is
``the'' totally definite algebra of
discriminant $1$ (or both). It is clear,
however, that matters can get complicated
when $F$ has many primes above $2$.}
\end{remark}

\section{Some special cases}
Let $X$ be as before a geometrically connected
component of $S\backslash\Shimp$. We assume
throughout this section that we are in the curve
case. Then for certain classes of examples we
will determine
the {\em deficient} local fields for $X$ in the
sense of \cite{JL1}, namely those local fields
$L$ containing the field of definition
of $X$ for which $X(L) = \emptyset$. In what
follows we shall make the following assumption

\begin{itemize}
\item[(Max)]\qquad $S$ is
a maximal compact subgroup of $\cG = \cG'$.
\end{itemize}

\noindent
All such $S$\/'s are conjugate in $\cG$.
In fact they are the units of
$\cM^\rint \otimes_{\cO_F}
    \bbO_\rf^{\cP}$
for any maximal $\cO_F$\/-order $\cM^\rint$ of
$B^\rint$. Likewise the maximal $S$\/'s are also
the units of
$\cM \otimes_{\cO_F}
    \bbO_\rf^{\cP}$
for any maximal $\cO_F[1/\cP]$\/-order $\cM$
of $B$. Such $S$\/'s certainly satisfy
Assumption (*).

The conjugacy classes of the $\cM^\rint\,$\/'s in
$B^\rint$ and of the $\cM$\/'s in $B$ are
classified by appropriate quotients of
$\Cl(F)/2\Cl(F)$: (See e. g. \cite[Section 3]{CF}
for the case of the $\cM^\rint\,$\/'s).
To avoid complications
we shall therefore add the following assumption:

\begin{itemize}
\item[(Odd)] \qquad The narrow class number $h^+(F)$
is odd.
\end{itemize}

\noindent
We list some consequences of Assumption~(Odd)
in the following
\begin{lemma}
\label{odd}
(a)\ For every $1\leq i\leq g(=[F:\bbQ])$ there
exists a unit $\epsilon_i$ of $F$ satisfying
$\infty_i(\epsilon_i) > 0$ and
$\infty_j(\epsilon_i) < 0$ for any $j\neq i$,
$1 \leq j \leq g$.

\noindent
(b)\ Totally positive units in $F$ are squares.

\noindent
(c)\ All the $\cM^\rint\,$\/'s are conjugate
in $B^\rint$ and all the $\cM$\/'s are
conjugate in $B$.

\noindent (d)\ All the geometrically connected
components of $S\backslash\Shimp$ are isomorphic.

\noindent
(e) The quotient of the normalizer of $\cM^\rint$
in $B_+^\rint$ by $Z(B^\rint)\m(\cM^\rint)\m$ is a finite group $(\bbZ/2\bbZ)^r$, where
$r$ is the number of primes dividing $\Disc B^\rint$, and by $B_+^\rint$ we denote elements of
$B^\rint$, whose norm is totally positive.
\end{lemma}
\begin{pf}
Though the proof is straightforward, we will sketch it for the convenience of the reader.
Our assumption (Odd) is equivalent to the decomposition
\begin{equation} \label{dec1}
\bbA\m=(\bbA\m)^2\bbO_{\rf}\m F\m.
\end{equation}
Then $-\epsilon_i$ is an element of
$F\m$ appearing in the decomposition of the idele $-1\in F_{\infty_i}\m\subset\bbA_F\m$, implying (a).
(b) now follows from (a) together with the fact that $\{-\epsilon_1,\ldots,-\epsilon_g\}$
form a basis of the $\bbF_2$-vector space $\cO_F\m/(\cO_F\m)^2$.

Decomposition (\ref{dec1}) implies decomposition $\bbA_{\rf}\m=(\bbA_{\rf}\m)^2\bbO_{\rf}\m F_+\m$.
Using strong approximation theorem together with (an analog of) (\ref{hasse}), we get that

\begin{equation} \label{dec2}
{\bf G^\rint}(\bbA_{F,\rf})={\bf Z(G^\rint)}(\bbA_{\rf}){\bf G^\rint}(\bbO_{\rf})
{\bf G^\rint}(F)_+,
\end{equation}

\noindent where ${\bf G^\rint}(F)_+=B_+^\rint$, and similarly for ${\bf G}$.
This decomposition immediately implies part (d) and reduces the remaining parts to the
corresponding local statemnts, which are clear.
\end{pf}

Let $X$ be a connected component of
$S\backslash\Shimp$. The field of definition $F'$
of $X$ is a Hilbert class field of $F$. Hence $F'$ is an abelian extension
of $F$ of odd degree. In particular, it is totally real.
$X$ has good reduction at all the (finite) primes of $F$ which do not divide
$\Disc B^\rint$. As in \cite{JL1}, the question of
existence of local points can be settled in such
cases via the trace formula and Hensel's lemma. In
\cite{Shi} Shimura had settled the case of real
points. His general results specialize in our case
to the following

\begin{proposition}
\label{real}
Let $\infty$ be a place of
$F'$ above $\infty_i$. Then
$X(F'_\infty) \neq \emptyset$ if and only if
$F(\sqrt{\epsilon_i}\,)$ does not split at
any prime where $B^\rint$ ramifies, where
$\epsilon_i$ is as in Lemma~\ref{odd}.
(Equivalently, if no finite prime of $F$ dividing
$\Disc B^\rint$ splits in $F(\sqrt{\epsilon_i}\,)$).
\end{proposition}

\begin{pf}
For the convenience of the reader and for
comparison with the case of a finite prime
$v|\Disc B$, we shall briefly sketch the proof.

The theorem on conjugation of Shimura varieties
allows us to assume that $B^\rint$ splits at
$\infty_i$. Putting $\cH = \bbC\smallsetminus \bbR$ we then
get, using Lemma~\ref{odd} (d), that
$(X\otimes_{F'} \bbC_\infty)^{{\rm an}} \simeq
\Gamma(S) \backslash \cH$. Then
\[ X(F'_\infty) \neq \emptyset
   \quad \Leftrightarrow \quad
   (\exists x\in \cH, \exists \gamma \in\Gamma(S):
    \gamma(x) = \ov{x}). \]

In this case, $\gamma^2(x) = x$, hence some power
$\gamma^k$ of $\gamma$ is central. Also
$\det\gamma < 0$, so $\gamma$ is conjugate in
$\GL_2(F'_\infty) = \GL_2(F_{\infty_i})$ to a
matrix $\matr{\alpha}{0}{0}{-\beta}$ with
$\alpha, \beta > 0$. Then $(\alpha/\beta)^k =1$
since $\gamma^k$ is central, so that $\alpha = \beta$.
It follows that $\epsilon:= \gamma^2 = \alpha^2$ is
in $\Gamma(S)\cap F = \cO_F^\times$. By
Lemma~\ref{odd} (b), we have
$F(\sqrt{\epsilon_i}\,) = F(\sqrt{\epsilon}\,)$,
so that this field splits $B^\rint$ and hence
cannot be split (over $F$) at any prime dividing
$\Disc B^\rint$.

Conversely, assume that $F(\sqrt{\epsilon_i})$ is not
split at any prime dividing $\Disc B^\rint$.
Then $F(\sqrt{\epsilon_i}\,)$ embeds into
$B^\rint$, and the image $\gamma$ of
$\sqrt{\epsilon_i}$  must belong to some maximal
order of $B^\rint$. By Lemma~\ref{odd} (c), we may
assume $\gamma\in\Gamma(S)$. Moreover in
$\GL_2(F_{\infty_i})$ $\gamma$ is conjugate to
a matrix $\matr{\alpha}{0}{0}{-\alpha}$, since
$\gamma^2 =\epsilon_i$ is a scalar negative at
$\infty_i$. Hence there exists $x\in \cH$
such that $\gamma(x) = \ov{x}$, therefore
$X(F'_\infty) \neq \emptyset$.
\end{pf}

Lastly,  Proposition~\ref{crit} gives the following
result concerning the deficiency of primes
$\cP|\Disc B^\rint$:

\begin{proposition}
\label{crit1}
Let $X$ be as before, and suppose that
Assumptions~(Odd) and (Max) hold and that
$fek_+$ is odd. Let $\cP$ be a finite prime
dividing $\Disc B^\rint$, and let $B$ be the
totally definite quaternion algebra over $F$
of discriminant $\Disc B^\rint = \Disc B/\cP$.
Let $\pi$ be a totally positive generator of
the principal ideal $\cP^{k_+}$.
Then in the notation of
Proposition~\ref{crit} we have that
$X(L)\neq\emptyset$ if and only if either
of Conditions (a) or (b) below holds:
\begin{enumerate}
\item[(a)] $F(\zeta_{2^{t(F)}})$ splits $B$, and
$\cP$ lies above $\tau$ with odd ramification
index. (Recall that $\tau$ is the prime above $2$
in $\bbQ(\zeta_{2^{t(F)}})^+$, where $t(F)$ is the
maximal integer for which $F$ contains
$\bbQ(\zeta_{2^{t(F)}})^+$.)
\item[(b)] $B$ is isomorphic to $B(-1,-\pi)$.
\end{enumerate}
\end{proposition}

\begin{pf}
We will check that our conditions (a) and (b) are equivalent to the corresponding conditions
of Proposition~\ref{crit}. The equivalence of conditions (b) follows from
Remark~\ref{quaternion} and Lemma~\ref{odd} (b). Since $t(F)\geq 2$, Lemma~\ref{rootof1} (ii)
implies that our condition
(a) follows from that of  Proposition~\ref{crit}.

Finally assume our assumption (a) and
choose an embedding of  $F(\zeta_{2^{t(F)}})$ into $B$. Observe that $\gm:=1+\zeta_{2^{t(F)}}$
is an algebraic integer such that $\Nm\gm\in\cO_{F,\cP}\m$.
Lemma~\ref{odd} (c) implies that some conjugate $\gm''$ of $\gm$ belongs to  $\GM(S)$.
By Lemma~\ref{rootof1}, $\gm''$ is of order $2^{t(F)}$ modulo the center, and by
Lemma~\ref{valp},
$$\val_\cP \Nm(\gm'')= e(\cP/\tau)\val_{\tau}(1+\zeta_{2^{t(F)}})(1+\ov{\zeta}_{2^{t(F)}})=e(\cP/\tau)$$
is odd. Hence using Lemma~\ref{kkplus} we can modify $\gm''$ by an element of $Z\GM(S)$
to get an element $\gm'$ satisfying $\val_\cP \Nm(\gm')=k$.
\end{pf}

\section{Evenness of jacobians of Shimura curves}

In \cite{PS} Poonen and Stoll defined a dichotomy
of principally polarized abelian varieties over a
number field $E$ into even and odd cases. When the
abelian variety is a jacobian of a curve $C/E$ they
gave a criterion for the evenness in terms of $C$.
Keeping assumptions (Max) and (Odd), let $X$ be a
component of an appropriate Shimura curve. The
methods of \cite{JL2} carry as are to our context.
We will call a prime of $F$ relevant if it is either infinite of divides
 $\Disc B^\rint$.
As in \cite[Theorem 23]{PS} it follows from
Corollary 10 there that the jacobian of $X$ is even
unless the genus of $X$ is even and the number of
relevant deficient primes for $X$ is
odd. We shall prove the following
\begin{theorem}
\label{PS} With $F$ quadratic, suppose that
assumptions (Max) and (Odd) hold. For technical
reasons assume also that $F$ is not isomorphic to
$\bbQ(\sqrt{2}\,)$ and that $\Disc B^\rint$ is not
a prime of residual characteristic $2$. Then the
jacobian of $X/F'$ is even.
\end{theorem}
\begin{pf}
We shall need the following
\begin{proposition}
\label{def}
 Let $F=\bbQ(\sqrt{m}\,)$ for a prime
$m\equiv 1\mod{4}$. Suppose $\Disc B^\rint$ is an
odd prime $\cP$ of $F$. Let $k_+$ be the order of
$\cP$ in $\Cl^+(F)$, and let $\pi$ be a totally
positive generator of  $\cP^{k_+}$.
Then\\
(i) The number of infinite deficient primes of $X$
is even if and only if $(-1,-\pi)_{F_\cP}= 1$.\\
(ii) The number of finite relevant deficient primes of $X$
is even if and only if $(-1,-\pi)_{F_\cQ}= 1$ for
$\cQ=\cP$ and for all $\cQ$ of residue
characteristic $2$.
\end{proposition}
\begin{pf}
As $F'$ is an abelian extension of $F$ of odd degree,
the number of primes $\infty'_1$, $\infty'_2$, and
$\cP'$ of $F'$ above each of $\infty_1$,
$\infty_2$, and $\cP$ is odd and they are deficient
simultaneously. For the infinite primes, the fields
$F(\sqrt{\epsilon_i}\,)$ are unramified at $\cP$
since $\cP$ is odd. By Proposition~\ref{real}, the
number of deficient infinite primes is even if and
only if $\cP$ is split in both or inert in both
these fields. This happens if and only if
$\left(\frac{-1}{\cP}\right)= \left(\frac{\epsilon_1}{\cP}\right)
\left(\frac{\epsilon_2}{\cP}\right)=1$, which
is equivalent to  $(-1,-\pi)_{F_\cP}= 1$, as $k_+$ and $\cP$ are odd.

For (ii) notice as before that only case (b) of
Proposition~\ref{crit1} can happen since $\cP$ is
odd. Thus the number of finite relevant deficient primes of
$X$ is even if and only if the quaternion algebra
$B(-1,-\pi)$ splits at all the finite places of
$F$.  The assertion now follows from
Remark~\ref{quaternion}.
\end{pf}
\begin{corollary}
(a) If $m \equiv 5\mod{8}$ then the number of
relevant deficient primes for $X$ is even.\\
(b) If $m \equiv 1\mod{8}$ and $\cP =p\cO_F$ with
$p$ a rational prime (inert in $F$, i. e.
$\left(\frac{m}{p}\right) = -1$) which is $\equiv
1\mod{4}$, then the number of relevant deficient primes for
$X$ is odd.
\end{corollary}
\label{def1}
\begin{pf}
(a) The condition on $m$ means that the rational
prime $2$ is inert in $F$. By the product formula,
$(-1,-\pi)_{F_\cP} = (-1,-\pi)_{F_2}$, so the
assertion follows.\\
(b) Since $p$ is inert in $F$ then without loss of
generality $\pi = p^{k_+}$\/. Also $-1$ is a square in
$F_p = F_\cP$, so $(-1,-\pi)_{F_\cP} = 1$. By the
proposition, the number of infinite deficient primes
is even. Now let $\lam$ be a prime of $F$ above
$2$. As $2$ splits in $F$ and $k_+$ is odd, we get
\[(-1,-\pi)_{F_\lam} = (-1,-p)_{\bbQ_2} = (-1,-p)_\bbR
(-1,-p)_{\bbQ_p} = -1,\]
\noindent and the claim follows from the Proposition.
\end{pf}

We now prove Theorem~\ref{PS}. The discriminant of
$B^\rint$ is a product of an odd number $r$ of
finite primes of $F$.
By Lemma~\ref{odd} (e), the group $W \simeq(\bbZ/2\bbZ)^r$
acts naturally on $X$. The stabilizer of a point
$x\in \cH$ in $W$ must be cyclic, hence trivial or
of order $2$. Let $n$ be the number of points of
$X$ fixed by a nontrivial element of $W$. Then the
genera $g(X)$ and $g(X/W)$ of $X$ and of $X/W$ are
related by
\[2-2g(X) = |W|(2-2g(X/W)) - n \;|W|/2.\]
Hence if $r\geq 3$ it follows that $4|2-2g(X)$ so
that $g(X)$ is odd and the jacobian $\Jac(X)$ is
certainly even. Now assume $r=1$ so that $\Disc
B^\rint = \cP$ for a prime $\cP$ of $F$ of residual
characteristic $p$. By genus theory (\cite[Ch. 3.8,
Cor.]{BSh}), Assumption (Odd) forces
$F=\bbQ(\sqrt{m}\,)$ for a prime $m$ which is either $2$ or
$\equiv 1\pmod{4}$. The case $m=2$ is excluded by our assumption. 
If the genus of $X$ is odd, there is nothing to prove. Else we must be in
either case 1.(a) or 1.(b)(i) of \cite[Theorem
1.1]{Sad}.

In the first of these cases $m\equiv 5\mod 12$ so
that $3$ is inert in $F$ and $\cP=3\cO_F$. Then
$(-1,-3)_\cP=1$, so that the number of deficient
infinite primes is even. We need therefore to show
that the number of relevant deficient finite primes is even
as well; equivalently, by Proposition~\ref{def},
that  $(-1,-3)_\cQ=1$ for all finite primes $\cQ$
of residue characteristic $2$. Set $(-1,-3)_2 =
\prod_\cQ (-1,-3)_\cQ$, where the product is over
$\cQ$ of residue characteristic $2$. By the product
formula $(-1,-3)_2=1$, so that if $2$ is inert we
are done. If $2$ is split (the only case left),
then for a prime $\cQ$ above $2$ we have $F_\cQ =
\bbQ_2$, so that
\[(-1,-3)_\cQ=(-1,-3)_2=
(-1)^{\frac{(-1-1)}{2}\frac{(-3-1)}{2}} = 1,\]
proving again what we need.

In the other case, $m\equiv 5\pmod{8}$, so
corollary~\ref{def1}(a) shows that the number of
relevant deficient primes is even, concluding the proof of
the Theorem.
\end{pf}


\end{document}